\documentclass{article}
\usepackage{amsfonts,amssymb,latexsym}
\newcounter{num}[section]
\newcommand{\Num}{\refstepcounter{num}%
	\textbf{\arabic{section}.\arabic{num}}}

\newcommand{\Theorem}{\textbf{Theorem~}}
\newcommand{\Proof}{{\textbf{Proof}}}
\newcommand{\Def}{\textbf{Definition~}}

\newcommand{\Lemma}{ \textbf{Lemma~}}

\newcommand{\Remark}{\textbf{Remark}}
\newcommand{\Prop}{\textbf{Proposition~}}
\newcommand{\Cor}{ \textbf{ Corollary~}}

\newcommand{\Ch}{{\mathfrak S}}

\newcommand{\Ac}{{\mathcal A}}
\newcommand{\Bc}{{\mathcal B}}
\newcommand{\Cc}{{\mathcal C}}
\newcommand{\Uc}{{\mathcal U}}
\newcommand{\Ec}{{\mathcal E}}
\newcommand{\Hc}{{\mathcal H}}
\newcommand{\Dc}{{\mathcal D}}
\newcommand{\Vc}{{\mathcal V}}
\newcommand{\Rc}{{\mathcal R}}
\newcommand{\Lc}{{\mathcal L}}
\newcommand{\Kc}{{\mathcal K}}
\newcommand{\Xc}{{\mathcal X}}
\newcommand{\Jc}{{\mathcal J}}
\newcommand{\Ic}{{\mathcal I}}

\newcommand{\al}{{\alpha}}
\newcommand{\la}{{\lambda}}
\newcommand{\La}{{\Lambda}}
	
\newcommand{\Fq}{\mathbb{F}_q}

\newcommand{\Gb}{\mathbb{G}}
\newcommand{\Lb}{\mathbb{L}}
\newcommand{\Ub}{\mathbb{U}}
\newcommand{\Ib}{\mathbb{I}}
\newcommand{\Jb}{\mathbb{J}}

\newcommand{\Hb}{\mathbb{H}}

\newcommand{\Ad}{{\mathrm{Ad}}}
\newcommand{\Mat}{{\mathrm{Mat}}}
\newcommand{\diag}{{\mathrm{diag}}}

\newcommand{\tx}{{\mathfrak t}}
\newcommand{\ux}{{\mathfrak u}}
\newcommand{\vx}{{\mathfrak{v}}}

\newcommand{\utx}{{\ux\tx}}

\newcommand{\eps}{{\varepsilon}}

\newcommand{\UT}{{\mathrm{UT}}}
\newcommand{\tUb}{\widetilde{\Ub}}
\newcommand{\GL}{{\mathrm{GL}}}
\newcommand{\Irr}{{\mathrm{Irr}}}
\newcommand{\Ind}{{\mathrm{Ind}}}

\newcommand{\row}{{\mathrm{row}}}
\newcommand{\col}{{\mathrm{col}}}

\newcommand{\Dp}{\Delta^+}

\renewcommand{\geq}{\geqslant}

\begin{document}
\Large

\title{Two supercharacter theories for the parabolic subgroups in  orthogonal and symplectic groups }
\author{A.N.Panov\\
\\
Department of Mathematics, 
Samara  University,\\
ul.~Akademika Pavlova 1, Samara, 443011, Russia\\
email: apanov@list.ru}
\date{}
 \maketitle
 \begin{abstract}
 	We construct two supercharacter  theories (in the sense of P. Diaconis and I.M. Isaacs) for the parabolic subgroups in  orthogonal and symplectic groups. For each supercharacter theory, we obtain a supercharacter analog of the A.A.Kirillov  formula  for irreducible characters of finite unipotent groups. 
 	
 \end{abstract}
{\small \textit{2000 Mathematics Subject Classification}. 20C33, 05E10. }

\section{Introduction. Main definitions and notations}

Traditionally, the following problem is considered as a main problem in the representation theory of finite groups:  to classify all irreducible representations  of a given finite group. As it turns out, for some groups (for example, the unitriangular  group) this problem is extremely  difficult.
 In the series of papers in  1995-2003  C.A.M. Andr\'{e} \cite{A1,A2,A3}  constructed the system of characters of the unitriangular group;
 he refered to these characters as the basic characters. Although these characters are reducible, they  have a range of properties
 similar to the ones of irreducible characters.
  Alongside, the system of basic varieties was constructed on which  the basic characters are constant  (the analog of conjugacy classes). 

  In  the paper  \cite{DI} (2008), P. Diaconis and I.M. Isaacs  introduced the  notion of a supercharacter theory  and constructed the supercharacter theory for algebra groups (by definition, an algebra group is the group of the form  $1+\Jc$, where $\Jc$ is an associative nilpotent finite dimensional algebra over a finite field).  The special case of this supercharacter theory is the theory of basic characters for the unitriangular group of C.A.M. Andr\'{e} in the interpretation of Ning Yan \cite{Yan}.
 Each group may afford different supercharacter theories. One can define a partial order on the set of supercharacter theories: a supercharacter theory is coarser than an  other one if its superclasses are unions of the superclasses of  other one. They set up a problem: given a group, to construct a supercharacter theory as close to the theory of irreducible characters as possible. Numerous papers are devoted to this problem   and also to solutions of representations theory problems in term of supercharacters  (see surveys  \cite{VERY, P0}).
 
   The supercharacter theory for Sylow subgroups of the orthogonal and symp\-lec\-tic groups was constructed in the papers  \cite{AN-1,AN-2,AFN,SA,P1}. This supercharacter theory is derived from the basic character  theory for the unitriangular group.

  The present paper is devoted to supercharacter theories of the parabolic subgroups in the orthogonal and symplectic groups. A parabolic subgroup is a semidirect product  $G=LU$ of a reductive group  $L$ and  unipotent normal subgroup  $U$.  The supercharacter theory for $U$ is constructed analogically to the Sylow subgroup.  It can  be assumed that there exists a supercharacter theory of $G$ constructed by the Mackey method: for the supercharacter  $\zeta$ of the subgroup  $U$, we consider its stabilizer  $H$ in $L$;  if  $\zeta$ can be extented to a character  of  $H\ltimes L$, then we construct the system of supercharacters  $\chi_{\theta,\zeta}=\Ind(\theta\cdot\zeta, H\ltimes L,G)$, where $\theta\in \Irr(H)$. However, these characters may not be orthogonal for different $\theta$. 
  For example, if  $G$ is a Borel subgroup in $\GL(4,\Fq)$, and  $\zeta $ is the supercharacter of  $U=\UT(4,\Fq)$ associated  with the linear form $\la=E_{13}^*+E_{24}^*$, then $H$ coincides with  $ \{\diag(a,b,a,b):~ a,b\in \Fq^*\}$ (stabilizer of  $\la$ in $L$);   easy to  prove that for $\theta\ne\theta'$ the characters   $\chi_{\theta,\zeta}$ and $ \chi_{\theta',\zeta}$ are not orthogonal if   $\theta$ and $\theta'$ are equal on the subgroup of scalar matrices. 
  Nevertheless, we  obtain an orthogonal system  of characters  replacing $H$ by  the smaller subgroup of scalar matrices. The last subgroup coincides with the intersection of all stabilizers of  $U\la U$.  This example gives rise an idea of  construction of a supercharacter theory by replacing the stabilizers of characters by their smaller subgroups.
  This idea  is implemented to construct  supercharacters theories for semidirect products with a normal algebra group  \cite{P4} and for parabolic subgroups in the present paper.  
  
  As a result,   we construct two  supercharacter theories  for the parabolic subgroups in the orthogonal and symplectic groups. These supercharacter theories are referred to as the $\Ub$-supercharacter theory and the $\Gb$-supercharacter theory.    In the first one, the supercharacters are construted by 
the $\Ub$-orbits in  $\ux^*$, in the second one -- by the  $\Gb$-orbits. Each  $\Gb$-orbit decomposes into the   $\Ub$-orbits; therefore, the  $\Gb$-supercharacter theory is coarser than the  $\Ub$- one.
 However, unlike the $\Ub$-supercharacter theory,  the superclasses in the $\Gb$-supercharacter theory afford  complete classification (see Propositions  \ref{Gequi} and \ref{Gstarequi}). 
  The solution of this classification problem in the  $\Ub$-supercharacter theory depends on  classification of elements in the Lie algebra of unipotent  radical of  parabolic subgroup for  $A$ series with respect to the following equivalence relation:  $x\sim x'$ if there exist the elements  $a,b\in \Ub$ and $r\in \Lb$ such that  $x'=ra x br^{-1}$.
   In the paper \cite{P3}, a conjecture  was proposed  and it was solved in the case of parabolic subgroups with blocks of sizes less or equal to 2.  
   
   The $\Ub$- and $\Gb$-supercharacter theories  provide the better approximation
    of the theory of irreducible characters than the general construction from the paper  \cite{H}.
   Observe that  the   $\Ub$- and  $\Gb$-orbits  in  $\ux^*$ are decomposed into the coadjoint orbits; according to the A.A.Kirillov orbit method, they parametrize the irreducible representations. As mentioned above, the classification of irredu\-cible representations and  coadjoint orbits for the unitriangular group is an extremely  difficult, "wild"\, problem.
   
  The main  results of the present paper are formulated in Theorems    \ref{supcharprime} and \ref{superchar}. We develop  the approach  for the parabolic subgroups of the general linear groups proposed by the author in the paper \cite{P4}. The constructed supercharacter theories are incomparable with the ones from the papers  \cite{P2,P3}. 
  
 The author suggests that the approach proposed in this paper could be applied to construct supercharacter theories for other groups such as   the contractions of simple groups of Lie type   \cite{Feigin, Panyu} and  the reductive groups over rings  \cite{Lu, St}.

 Let us give the definition of  a supercharacter theory (see \cite{DI}).
Let  $G$ be a finite group, $1\in G$ be the unit element. Let  $\Ch = \{\chi_1, \ldots, \chi_M\}$ be a system of characters (representations) of the group  $G$. \\
 \Def\Num\label{superch}. The system of characters $\Ch$  determines a supercharacter theory of $G$ if there exists a partition $\Kc = \{K_1,\ldots, K_M\}$ of the group  $G$ satisfying the following conditions: \\
 S1) the characters of $\Ch$ are pairwise disjoint  (orthogonal);\\
 S2) ~ each character  $\chi_i$ is constant on each subset $K_j$;\\
 S3) ~ $\{1\} \in \Kc$.
 \\
 Each character of $\Ch$  is  referred to as  a  {\it  supercharacter}, each subset of $\Kc$  --  a {\it  superclass}. Observe that the number of supercharacters is equal to the number of superclasses.

 For each supercharacter  $\chi_i$, consider the subset  $X_i$ of all its irreducible constituents.  Observe that the condition S3)  of Definition  \ref{superch} may be replaced by  following condition:\\
 S3')  The system of subsets  $\Xc=\{X_1,\ldots, X_N\}$ is a partition  of the system of irreducible characters  $\Irr(G)$. Moreover,
 here each character   $\chi_i$ differs from the character  $\sigma_i=\sum_{\psi\in X_i} \psi(1)\psi$ by a constant factor (see \cite{DI}).

Let  $\Fq$ be a finite field of $q$ elements and characteristic  $>2$.  Denote by $\Ib_N $ the $N\times N$ matrix in which all entries on the anti-diagonal are equal to $1$ and all other entries are zeros. Consider the involutive antiautomorphism $X^\dag= \Ib_{N} X^t \Ib_{N}$ of the algebra of  $N\times N$ matrices over the field  $\Fq$. The orthogonal group of type $B_n$  (respectively, $D_n$) consists of matrices   $g\in\GL(N,\Fq)$ obeying  $g^\dag=g^{-1}$, where $N=2n+1$ (respectively,  $N=2n$).   The symplectic group of type  $C_n$ is defined similarly, applying the involutive antiautomorphism   $X^\dag= \Jb_{2n}X^t\Jb_{2n}$, where $N=2n$ and $\Jb_{2n}=\left(\begin{array}{cc} 0&\Ib_n\\ -\Ib_n&0\end{array}\right).$
 
 We enumerate the rows and columns of  $N\times N$  matrices as follows:
 $$
 \begin{array}{c}
  n>\ldots > 1 > 0 > -1> \ldots >-n\quad \mbox{for} \quad N=2n+1,\\ n>\ldots > 1> -1> \ldots >-n \quad \mbox{for} \quad
 N=2n.
 \end{array}
 $$ 
 
 A parabolic subgroup  $\Gb$  in   $\GL(N,\Fq)$ is defined by the decomposition  $\Ic$ of segment  $I=[-n,n]$ in the case  $B_n$ (respectively, $I=[-n,n]\setminus \{0\}$ in the cases $C_n$ and $D_n$) into a union of consecutive segments.

 We enumerate the components of decomposition  
  \begin{equation}\label{decomI}
 \Ic=\{I_{\ell},\ldots, I_0,\ldots , I_{-\ell}\}
 \end{equation} 
  comparably with the enumeration of rows amd columns. The subgroup $\Gb$ is a semidirect product  $\Gb=\Lb\Ub$ of Levi subgroup   $\Lb =\GL(n_\ell)\times\cdots\times \GL(n_0,\Fq))\times\cdots\times  \GL(n_{-\ell})$, where $n_i=|I_i|$, and the radical $\Ub$. The subgroup  $\Ub$ is an algebra group $\Ub=1+\Uc$, where $\Uc$ is a nilpotent associative subalgebra in $\Mat(n,\Fq)$. The subgroup  $\Ub$ acts on $\Uc$ by the left and right multiplication. These actions induces the action of  $\Ub$ on the dual space  $\Uc^*$ by  the formulas  
 $a\La(x)=\La(xa)$ and $\La a(x)=\La(ax)$, where $\La\in \Uc^*$,~ $x\in\Uc$,~  $a\in\Ub$.
 
 Suppose that the parabolic subgroup  $\Gb$ is invariant with respect to  $A\to A^\dag$. This is equivalent to the condition that the decomposition  $\Ic$ is symmetric relative to zero.  Then the subgroup   $G$ of matrices  $g\in\Gb$,~ $g^\dag=g^{-1}$,  is a parabolic subgroup in the orthogonal (symplectic ) group.
 The subgroup $G$ is a semidirect product $G=LU$, where $L$ and $U$ are intersections of  $\Lb$ and $\Ub$ with the orthogonal (symplectic) group.   
  The subgroup  $L$ consists of matrices $\mathrm{diag}(A_{\ell},\ldots, A_0, \ldots, A_{-\ell})$, where   for each $\ell\geq i\geq {-\ell}$ the transpose of the matrix   $A_i$ about its anti-diagonal coincides with  $A_{-i}^{-1}$.  
The subgroup $U$ is an  unipotent subgroup with the Lie algebra  $\ux=\{x\in \Uc: ~~ x^\dag =-x\}$.

 The group $\Gb$ (respectively, the subgroup $\Ub$) acts on  $\ux$ by the  formula $g\centerdot x= gxg^\dag$, where $g\in \Gb$ and $x\in\ux$. 
Restricting  this action on  $\Ub$, we obtain the action of the subgroup  $\Ub$ on $\ux$.
Observe that if  $u\in U$, then  $u\centerdot x$ coincides with the adjoint action $\Ad_u(x)$. 

The action of $\Gb$ provides the action of  $\Gb$ on $\ux^*$ by the formula  $g\centerdot \la(x)=\la(g^\dag x g)$. Similarly, one can define the action of  $\Ub$ on $\ux^*$. 

To construct supercharacter theory we introduce the notion of a Springer map.\\
\Def\Num\label{Spr} (\cite{SA}). 
A map  $f:\Ub\to\Uc$ is called a Springer map if $f$ is a bijection obeying the following conditions:\\
1)  ~ $f(U)= \ux $,\\
2) ~ there exist  $a_2,a_3,\ldots $ from the field $\Fq$ such that  $f(1+x) = x+a_2x^2+a_3x^3+\ldots$ for any  $1+x\in \Ub$.

It follows that a Springer map satisfies the condition  $\Ad_vf(u)=f(vuv^{-1})$ for   $u,v\in \Ub$. 
Examples of a Springer map are as follows:\\
1) the logarithm map  $\ln(1+x) = \sum_{i=1}^\infty (-1)^{i+1}\frac{x^i}{i}$~~ (it requires strong restrictions on the  characteristic of  field);\\
2) the Cayley map $f(1+x)=\frac{2x}{x+2}$ ~~(for $\mathrm{char}\,\Fq\ne 2)$.

In what follows we fix a  Springer map on  $\Ub$. 

\section{ $\Ub$-supercharacter theory}

In this section, we use $\Ub$-orbits in  $\ux^*$ to construct a supercharacter theory for the parabolic subgroup  $G$. 
We apply the approach  that is used in the  paper  \cite{SA} for construction a supercharacter theory for Sylow subgroups in the  orthogonal and symplectic groups. 

let $\Hc$ be the subspace in   $\Uc$ that consists of all matrices  $x=(x_{ij})$ obeying the property: if $x_{ij}\ne 0$, then $i\in I_\ell\sqcup\ldots\sqcup I_0$.
The subspace  $\Hc$ is an ideal in the associative algebra $\Uc$. Respectively, $\Hb=1+\Hc$ is an algebra group and a normal subgroup in  $\Ub$.
For $\La\in\Uc^*$, the following associative subalgebras in $\Uc$ are defined:
$$\Rc_\La=\{x\in\Uc:~~ \La(x\Hc)=0\},\quad
\Lc_\La=\{x\in\Uc:~~\La(\Hc^\dag x)=0\},\quad
\Uc_\La=\Rc_ \La\cap\Lc_\La.$$
Then   $\Ub_\La=1+\Uc_\La$ is an algebra subgroup in  $\Ub$.
The restriction of Springer map on  $\Ub_\La$ defines  a bijection $f: \Ub_\La\to\Uc_\La$.
\\
\Lemma\Num\label{xy} (see \cite[Lemma 6.3]{SA}).  If  $x,y\in\Uc_\La$, then $\La(xy)=0$.

Let  $\Lb_\La$ be a stabilizer of  $\La\in\Uc^*$ with respect the coadjoint representation. Since  $\Hc$ and $\Hc^\dag$ are ideals in $\Uc$, we see $\Ad_r(\Uc_\La)=\Uc_\La$ for any  $r\in \Lb_\La$. 

For any subgroup  $\Lb_0$ in $\Lb_\La$, we construct  a semidirect product 
$\Gb_{0,\La}=\Lb_0\Ub_\La$.  
Let $T_r$, ~$r\in L$, be a representation of the subgroup  $\Lb_0$.
Lemma  \ref{xy} implies that  the formula
\begin{equation}\label{Xieps}
\Xi_{T,\La}(g)= T(r)\eps^{\La(f(u))},
\end{equation}
where $g=ru$,~ $r\in \Lb$, ~$u\in \Ub$,  defines a representation of the group  $\Gb_{0,\La}$.
If $\theta$ is the character of representation  $T$, then the formula  
\begin{equation}\label{xieps}
\xi_{\theta,\La}(g)= \theta(r)\eps^{\La(f(u))}.
\end{equation}
defined a character of the group   $\Gb_{0,\La}$.

Let $\la\in\ux^*$. There exists a unique element $\La\in\Uc^*$,~ $\La^\dag=-\La$, such that its restriction on  $\ux$ coincides with  $\la$.
Consider the subalgebra  $\ux_\la=\ux\cap \Uc_\La$. 
Easy to see   $\La(x\Hc) = \La(\Hc^\dag x)$ for any  $x\in \ux$. Therefore, 
$$ \ux_\la=\ux\cap \Rc_\La=\ux\cap \Lc_\La.$$
The Springer map  $f$ bijectively maps  $\Ub\to \Uc$ and $\Ub_\La\to \Uc_\La$; then  $f$ bijectively maps the subgroup  $U_\la=U\cap\Ub_\La$ onto $\ux_\la$. 

Let  $\Pi$ define the natural projection  $\Uc^*$ on $\ux^*$.\\
\Lemma\Num\label{Pi} (see \cite[Lemma 6.8]{SA}). $\Pi(\Hb \La)=\Hb\centerdot\la$.\\
\Proof.  For any  $h=1+y\in\Hb$, we obtain  $h\centerdot \la (x)= h\la h^\dag(x) = 
\la(h^\dag xh)=\la((1+y^\dag)x(1+y))=\la(x)+\La(y^\dag x) +\La(xy) + \La(y^\dag x y)$.

One can directly calculate that  $y^\dag x y=0$ for any  $y\in \Hc$ and $x\in\ux$.
As $\La^\dag=-\La$, for any  $x\in \ux$, we have   $$ \La(y^\dag x) = -\La^\dag(y^\dag x) = -\La((y^\dag x)^\dag)=-\La(x^\dag y)=\La(xy).$$
Hence  $h\centerdot \la (x)=\la(x)+2\La(xy)=(1+2y)\La(x)$. $\Box$

Let  $\pi$ define the natural projection  $\ux^*\to\ux_\la$.\\
\Lemma\Num\label{pi} (see the proof of \cite[Theorem 6.9]{SA}).  $\pi^{-1}(\la)=\Hb\centerdot \la$.\\
\Proof. Let  $P$  be the natural projection  $\Uc^*\to\Rc_\La^*$. 
 The definition of  $\Hc$ implies  $\Rc_\La^\perp=\Hc\La$.
Then $P^{-1}P(\La) = \La+\Rc_\La^\perp = (1+\Hc)\La=\Hb\La$. By Lemma \ref{Pi}, we get  $$\pi^{-1}(\la)=\Pi(P^{-1}P(\La))=
\Pi(\Hb\La) = \Hb\centerdot\la.~~\Box $$ 

Denote 
$$L_{\Ub\La\Ub}=\{h\in L:~~ \Ad^*_h(M)=M~~\mbox{for~any}~~ M\in\Ub\La\Ub\}.$$
Let $L_0$ be a subgroup in  $L_{\Ub\La\Ub}$.
For any character  $\theta$  of the subgroup  $L_0$, 
the formula  (\ref{xieps}) defines the character 
$\xi_{\theta,\la}$ of subgroup $G_{0,\la}=L_{0}U_\la$.
Consider the character  $\chi_{\theta,\la}$ induced from the character  $\xi_{\theta,\la}$ of  subgroup  $G_{0,\la}$ to the group  $G$.
  
 \Prop\Num\label{propchial}. The character  $\chi_{\theta,\la}$ is calculated by the formula
  \begin{equation}\label{chial}
  \chi_{\theta,\la}(g)=\frac{|\Hb\centerdot\la|}{|\Ub\centerdot\la|\cdot |L_0|}\cdot\sum_{\rho\in L} \dot{\theta}(\rho r \rho^{-1}) \sum_{\mu\in \Ub\centerdot \la} \eps^{\mu(f(\rho u\rho^{-1}))},
  \end{equation}
  where $g=ru,~ r\in L,~ u\in U$.\\
  \Proof. Consider the subgroup  $G_0=L_0U$.
  	The character  $\chi_{\theta,\la}$ is induced from the character  $$\chi^\circ_{\theta,\la} =
  	\Ind(\xi_{\theta,\la},G_{0,\la},G_0)$$
  	of  subgroup  $G_0$.
The value of the character  $\chi^\circ_{\theta,\la}$ on the element   $g_0\in G_0$ is calculated by the formula 
 $$\chi^\circ_{\theta,\la} (g_0)=\frac{1}{|G_{0,\la}|}\sum_{s_0\in G_0}\dot{\xi}_{\theta,\la}(s_0g_0s_0^{-1}),$$ 
where $\dot{\xi}_{\theta,\la}$ is the function on  $G_0$ that equals to  $\xi_{\theta,\la}$ on  $G_{0,\la}$  and equals to zero outside of  $G_{0,\la}$. 
The function   $\dot{\xi}_{\theta,\la}$ decomposes into a product of two functions   $\dot{\xi}_{\theta,\la} = \theta\dot{\eps}^{\la\cdot f}$, where $ \dot{\eps}^{\la\cdot f}$ is the function on  $U$ that equals to $\eps^{\la\cdot f}$ on $U_\la$ and equals to zero outside of  $U_\la$.

For $s_0=\rho_0v$ and $g_0=r_0u\in G_0$, where  ~$\rho_0, r_0  \in L_0$,~ $v, u\in U$, we obtain  
 $s_0g_0s_0^{-1} = (\rho_0 r_0 \rho_0^{-1})(\rho_0 v^{r_0}uv^{-1}\rho_0^{-1})$.
 We get
 \begin{equation}\label{circ}
 \chi^\circ_{\theta,\la} (g_0)=\frac{1}{|G_{0,\la}|}\sum_{\rho_0\in L_0,~v\in U}\theta(\rho_0 r_0 \rho_0^{-1}) 
 \dot{\eps}^{\la(f(\rho_0 v^{r_0}uv^{-1}\rho_0^{-1}))}.
 \end{equation}

By $\rho_0\in L_{\La}$, it follows 
$\la(f( \rho_0v^{r_0}uv^{-1}\rho_0^{-1})) = \la(f( v^{r_0}uv^{-1})).$

Let $v=1+y$,~ $y\in \Uc$. 
Since $r_0\in L_0\subseteq L_{\Ub\La\Ub}$, we verify $\La(w y^{r_0}w') = \La(w y w')$ for any matrices  $w$ and $w'$ from $\Ub$ or $\Uc$. Hence
\begin{equation}\label{rhozero}
\la(f( \rho_0v^{r_0}uv^{-1}\rho_0^{-1})) = \la(f( v^{r_0}uv^{-1})) = \la(f( vuv^{-1})).
\end{equation}

Substituting  (\ref{rhozero}) in (\ref{circ}), we have 
$$ \chi^\circ_{\theta,\la} (g_0) = \frac{|L_0|}{|G_{0,\la}|}\cdot\theta(r_0)\sum_{v\in U}
\dot{\eps}^{\la(f( vuv^{-1}))}.$$

It follows from Lemma \ref{pi} that 
$$\dot{\eps}^{\la(f(vuv^{-1}))}=\frac{|U_\la|}{|U|}\cdot\sum_{\mu\in \pi^{-1}(\la)}\eps^{\mu(f(vuv^{-1}))} =$$ 
$$\frac{|U_\la|\cdot |\Hb\centerdot \la|}{|U|\cdot |\Hb|}\cdot\sum_{h\in\Hb}\eps^{h\centerdot\la (f(vuv^{-1}))}=
\frac{|U_\la|\cdot |\Hb\centerdot \la|}{|U|\cdot |\Hb|}\cdot\sum_{h\in\Hb}\eps^{v^{-1}h\centerdot\la (f(u))}.$$
 
 We obtain
$$ \chi^\circ_{\theta,\la} (g_0) = \frac{|L_0|}{|G_{0,\la}|}\cdot \frac{|U_\la|\cdot |\Hb\centerdot \la|}{|U|\cdot |\Hb|}\cdot\theta(r_0)\sum_{v,\in U,~ h\in\Hb}\eps^{v^{-1}h\centerdot\la (f(u))} =$$
$$  \frac{ |\Hb\centerdot \la|}{|U|\cdot |\Hb|}\cdot\theta(r_0)\sum_{v,\in U,~ h\in\Hb}\eps^{v^{-1}h\centerdot\la (f(u))}$$

Easy to see that  $\Ub=U\Hb$. Then  
$$\chi^\circ_{\theta,\la} (g_0)=\frac{ |\Hb\centerdot \la|}{|U|\cdot |\Hb|}\cdot |U\cap \Hb|\cdot\theta(r_0)\sum_{w\in\Ub}\eps^{w\centerdot\la (f(u))}=$$
$$\frac{ |\Hb\centerdot \la|}{|U|\cdot |\Hb|}\cdot |U\cap \Hb|\cdot\frac{|\Ub|}{|\Ub\centerdot\la|}\cdot\theta(r_0)\sum_{\mu\in\Ub\centerdot\la}\eps^{\mu (f(u))}= 
\frac{ |\Hb\centerdot \la|}{|\Ub\centerdot\la|}\cdot\theta(r_0)\sum_{\mu\in\Ub\centerdot\la}\eps^{\mu (f(u))}.$$
This implies  formula  (\ref{chial}). $\Box$

Letting $L_0=\{1\}$ in the proof of Proposition  \ref{propchial}, we obtain the statement on the characters of subgroup  $U$. 
\\
\Prop\Num\label{uchi} (see \cite[Theorem 6.9]{SA}). Let $\zeta_\la$ be a character of subgroup  $U$ induced from the character  $\eps^{\la\cdot f}$ of subgroup $U_\la$. Then
 $$\zeta_\la(u)= \frac{|\Hb\centerdot\la|}{|\Ub\centerdot\la|}\sum_{\mu\in\Ub\centerdot\la}\eps^{\mu(f(u))}.$$
 The system of characters  $\{\zeta_\la\}$ and the partition $\{f^{-1}(\Ub\centerdot x)\}$, where $\la$ and  $x$ run through the set of representatives of  $\Ub$-orbits in  $\ux^*$ and $\ux$, respectively, give rise to a supercharacter theory for the group  $U$. $\Box$

Introduce notations\\
 $S_{\Ub\centerdot\la }=\{h\in L:~~ \Ad^*_h(\Ub\centerdot\la)=\Ub\centerdot\la\}$,\\
  $S_{\Ub\La\Ub}=\{h\in L:~~ \Ad^*_h(\Ub\La\Ub)=\Ub\La\Ub\}$.\\
 \Lemma\Num\label{ss}.  $S_{\Ub\centerdot\la}=S_{\Ub\La\Ub}$.\\
 \Proof. If $h\in S_{\Ub\centerdot\la}$, then $\Ad^*_h(\la)=a\la a^\dag$ for some  $a\in \Ub$. Hence $\Ad^*_h(\La)=a\La a^\dag$  and $h\in S_{\Ub\La\Ub}$.
On the other hand,  if $h\in  S_{\Ub\La\Ub}$, then $\Ad^*_h(\La)\in S_{\Ub\La\Ub}$.
 The element $\Ad^*_h(\La)$, as $\La$, is an invariant of the transformation  $M\to -M^\dag$ on the $\Ub-\Ub$ orbit of $\La$.  It follows from \cite[Cor. 13.9]{Isaacs} that the set of invariants consists of   $\{a\La a^\dag\}$. Then  $\Ad^*_h(\La) =a\La a^\dag$  and $\Ad^*_h(\la) =a\la a^\dag$. We have  $h\in S_{\Ub\centerdot\la}$. $\Box$\\
\Cor\Num.  $L_{\Ub\La\Ub}$ is a normal subgroup in  $S_{\Ub\centerdot\la}$

Denote by  $\Ac'_0$ the set of pairs  $\{(\theta,\la)\}$, where $\la$ runs through the set of representatives of  $\Ub$-orbits in  $\ux^*$, and  $\theta$ runs through the set of  $S_{\Ub\centerdot\la}$-irreducible characters of subgroup  $L_{\Ub\La\Ub}$. The group  $L$ acts on  $\Ac'_0$ by conjugation. Denote the set of $L$-orbits by  $\Ac'=\Ac'_0/L$.
We attach to  $\al=(\theta,\la)\in\Ac'$ the character defined by the formula  (\ref{chial}) for $L_0=L_{\Ub\La\Ub}$. We get  \begin{equation}\label{alchial}
 \chi_\al(g)=\frac{|\Hb\centerdot\la|}{|\Ub\centerdot\la|\cdot |L_{\Ub\La\Ub}|}\cdot\sum_{\rho\in L} \dot{\theta}(\rho r \rho^{-1}) \sum_{\mu\in \Ub\centerdot \la} \eps^{\mu(f(\rho u\rho^{-1}))}.
\end{equation} 
  We show  in Theorem   \ref{supcharprime} that   $\{\chi_\al:~ \al\in \Ac'\}$ is a system of supercharacters.

Our next goal is to construct  superclasses.
For  $h\in L$, we consider the   largest  $\Ub-\Ub$ invariant subspace   $\Vc^*_h$  in $\Uc^*$ on which  $\Ad_h^*$ acts identically.
Then its orthogonal complement  $\Uc_h=(\Vc^*_h)^\perp$ is the smallest  
$\Ub-\Ub$ invariant subspace in  $\Uc$ such that  $\Ad_h$ identically acts on $\Uc/\Uc_h$.
 The intersection  $\ux_h=\ux\cap\Uc_h$ is invariant with respect to the  $\Ub$-action on $\ux$,  and $\Ad_h$ is identical on $\ux/\ux_h$. 

Consider the set   $\Bc'_0$ of pairs  $\{(h,\omega)\}$, where  $h\in L$, and $\omega$ runs through the set of   $\Ub$-orbits in  $\ux/\ux_h$.  
 The group  $L$ acts on $\Bc'_0$ by conjugation. Denote  $\Bc'=\Bc'_0/L$.
 
We attach to  $\beta = (h,\omega)$ the class  
 $$K_\beta =\bigcup_{r\in L} r(h\cdot f^{-1}(\omega+\ux_h))r^{-1}.$$

\Theorem\Num\label{supcharprime}.  The system of characters  $\{\chi_\al:~\al\in\Ac'\}$  and the partition of parabolic subgroup  $G$ into classes  $\{K_\beta:~\beta\in \Bc'\}$ give rise to a supercharacter theory.  
\\
\Proof. 
\textit{Item 1.} Let us show  $|\Ac'|=|\Bc'|$.
Consider the set  $\Cc'_0$ of pairs  $(h,\la)$, where  $\la$ runs through the set of representatives of  $\Ub$-orbits in $\ux^*$, and $h\in L_{\Ub\La\Ub}$. As above 
$\Cc'=\Cc'_0/L$.

Then $|\Cc'|$ equals to  the sum over all $\Ub$-orbits in  $\ux^*$ of
the numbers of all $S_{\Ub\centerdot \la}$ orbits in $ L_{\Ub\La\Ub}$.
 This sum coincides with the sum over all  $\Ub$-orbits in $\ux^*$ of the numbers of all $S_{\Ub\centerdot\la}$-irreducible characters of   $L_{\Ub\La\Ub}$.  Then  $|\Cc'|=|\Ac'|$.

Let us show  $|\Cc'|=|\Bc'|$. For each  $h\in L$, we denote  $$\mathrm{Cent}(h)=\{r\in L:~ hrh^{-1}=r\}.$$ 
The number $|\Cc'|$ coincides with the sum over   $h\in L$ of the numbers of all $\mathrm{Cent}(h)$-orbits on the set of  $\Ub$-orbits such that $h\in L_{\Ub\La\Ub}$.
The condition  $h\in L_{\Ub\La\Ub}$ is equivalent to    $\La\in \Vc^*_h=\Uc_h^\perp$, which is equivalent to 
$\la\in \ux_h^\perp$. So, the number $|\Cc'|$ coincide with the sum over  $h\in L$ of the numbers of  all $\mathrm{Cent}(h)$-orbits on the set of  $\Ub$-orbits on $ \ux_h^\perp$.
In its turn, the number of  $\mathrm{Cent}(h)$-orbits on the set of  $\Ub$-orbits on  $ \ux_h^\perp$  coincides with the number of orbits of the semidirect product  
 $\mathrm{Cent}(h)\ltimes \Ub$ on $(\ux/\ux_h)^*$.   The orbit number  of any finite transformation  group of  a linear space over a finite field equals to the orbit number  in the dual space.  Therefore, the number  $|\Cc'|$ coincides with the sum over  $h\in L$  of the numbers of all
$\mathrm{Cent}(h)$-orbits on the set of  $\Ub$-orbits  in $\ux/\ux_h$. Hence  $|\Cc'|=|\Bc'|$. We obtain  $|\Ac'|=|\Bc'|$.\\
\textit{Item 2.} Let us show that the characters $\{\chi_\al\}$ are pairwise orthogonal. Since  $\{\eps^\mu\}$ is the system of characters on  abelian group  $\ux$,  these functions are pairwise orthogonal on  $\ux$. The Springer map  $f:U\to \ux $ is bijection; the functions   $\{\eps^{\mu\cdot f}\}$ are pairwise orthogonal on  $U$.
The functions of the form  $\sum_{\mu\in \Ub\centerdot\la}\eps^{\mu\cdot f}$ are orthogonal on $U$ if they are  attached to different  $\Ub$-orbits. 

Suppose that two different elements  $\al_1,\al_2$ from $\Ac'$ are presented by the pairs  $(\theta_1,\la_1)$ and  $(\theta_2,\la_2)$ from $\Ac'_0$.  If $\la_1$ and $\la_2$ belong to different  $L\ltimes \Ub$ orbits, then the characters  $\chi_{\al_1}$ and $\chi_{\al_2}$
are orthogonal because their restrictions on  $U$ are orthogonal.

If  $\la_1$ and  $\la_2$ belong to a common  $L\ltimes \Ub$ orbit, then it is sufficient to consider the case when they belong to a common  $\Ub$-orbit $\Ub\centerdot\la$. The characters $\chi_{\al_1}$ and $\chi_{\al_2}$ are orthogonal since the $S_{\Ub\centerdot\la}$-irreducible characters  $\theta_1$ and $\theta_2$ are orthogonal.\\
\textit{Item 3.} Let us show that the characters  $\{\chi_\al\}$  are constant on the classes
$\{K_\beta\}$. 
Consider the function 
$$F_{\theta,\la}(g)=\dot{\theta}(r)\sum_{\mu\in\Ub\centerdot\la}\eps^{\mu(f(u))}=
\dot{\theta}(r)\zeta_\la(u),$$
where $g=ru$,~$r\in L$,~$u\in U$,
and the subset $\kappa_{h,\omega}=h\cdot f^{-1}(\omega+\ux_h)$.
Since the character  $\chi_\al$ and the class $K_\beta $ are averages with respect to $\Ad_\rho$, ~ $\rho\in L$, of the function  $F_{\theta,\la}$ and the subset 
$\kappa_{h,\omega}$, respectively, it is sufficient to show that the function  $F_{\theta,\la}$ 
is constant on $\kappa_{h,\omega}$.
	
	If $h\notin L_{\Ub\La\Ub}$, then $F_{\theta,\la}(\kappa_{h,\omega})=0$.
Suppose  $h\in L_{\Ub\La\Ub}$. Then, as above,  $\la\in \ux_h^\perp$.
Hence $\mu(\omega+\ux_h^\perp)= \mu(\omega)$. We get
$$F_{\theta,\la}(\kappa_{h,\omega})=\theta(h)\zeta_\la(f^{-1}(\omega))=\mathrm{const} \Box.$$	
\Cor\Num. 1) Each irreducible character of $G$ is a constituent of a unique supercharacter   from  $\{\chi_\al:~\al\in\Ac'\}$.
2) Each conjugacy class of $G$ is contained in a unique class from    $\{K_\beta:~\beta\in \Bc'\}$.

\section{$\Gb$-supercharacter theory}

In this section, we use $\Gb$-orbits in  $\ux^*$ to construct a supercharacter theory for the parabolic subgroup  $G$.
The classification of $\Gb$-orbits  in $\ux$  and $\ux^*$ may be described  in terms of the basic subsets (rook placements) in the set of positive roots.

The set of matrix units  $E_{i,j}$, where $n\geq i>j\geq -n$, forms a basis in  $\utx(N,\Fq)$. 
We say that  $i=\row(\gamma)$ is the row number of the pair $\gamma=(i,j)$, and  $j=\col(\gamma)$ is the column number.

We refer to  $\gamma=(i,j)$, where  $n\geq i>j\geq -n$, as a \textit{positive root} if\\
$i>j>-i$ in the case $B_n$,\\
$i>j\geq -i$,~ $j\ne 0$ in the case  $C_n$,\\
$i>j> -i$,~ $j\ne 0$ in the case  $D_n$.\\
We denote the set of positive roots  by $\Dp$. 

For positive root  $\gamma=(i,j)$, we denote  $\gamma'=(-j,-i)$ (by definition, the pair  $\gamma'$ is not a positive root). Attach to  $\gamma\in\Dp$ the matrix $\Ec_\gamma=E_\gamma+\epsilon(\gamma)E_{\gamma'}$, where $\epsilon(\gamma)\in\{-1,0,1\}$, in the Lie algebra of Sylow subgroup.

Denote by  $\Dp_\ux$ the set of $\gamma\in \Dp$ for which $\Ec_\gamma\in \ux$.
The system of matrices   $\{\Ec_\gamma:~\gamma\in\Dp_\ux\}$  is a basis in  $\ux$.

For any subset  $\Dc\subset \Dp_\ux$ denote   $\Dc'=\{\gamma':~ \gamma\in\Dp_\ux\}$.\\
\Def\Num. We refer to  $\Dc$ as a basic subset  (rook placement) if there is no more than one pair of  $\Dc\cup\Dc'$ in each row and each column  (i.e.   $\Dc\cup\Dc'$  is a basic subset in the sense of C.A.M. Andr\'{e} in the set positive roots  $\{(i,j): n\geq i>j\geq -n\}$ for the unitriangular group   \cite{A1,A2,A3}).

Attach to each basic subset   $\Dc$  and a map $\phi:\Dc\to \Fq^*$  the element in  $\ux$ as follows 
\begin{equation}\label{dphi}
x_{\Dc,\phi}=\sum_{\gamma\in\Dc}\phi(\gamma)\Ec_\gamma.
\end{equation}

Fix the non-square element $\delta\in \Fq^*$. Observe that each root $\gamma\in\Delta^+_\ux$ belongs to a unique  $I_k\times I_m$  of components of decomposition (\ref{decomI}). We refer to $k$ as a \textit{block-row} (respectively, $m$ as a  \textit{block-column}) of the root $\gamma$.\\
\Def\Num\label{special}. We call the map $\phi$ special if \\
1)~ $\phi(\gamma)\in\{1,\delta\}$, \\
 2) ~  $\phi(\gamma)=\delta$  implies $\gamma=(i,-i)$ for some $n\geq i\geq 1$ (in the case  $C_n$),\\
 3) for  any  $\ell\geq k \geq 1$ there is no more than one root $\gamma= (i,-i)$ in $\Dc\cap\left(I_k\times I_{-k}\right)$ that obeys 2).
 
 We refer to a pair $(\Dc,\phi)$, where $\Dc$ is a basic subset and $\phi$ is special, as a \textit{basic pair.}
 
For special $\phi$ and  each  $\ell\geq k \geq 1$, we define the number $d_k(\phi)$ as follows: if the there exists  a root  $\gamma=(i,-i)\in I_k\times I_{-k}$ with $\phi(\gamma) = \delta$, we take $d_k(\phi)=-1$; if there is no such root, we take $d_k(\phi)=1$.  So in cases $B_n$ and $D_n$, we have $\phi(\gamma)=1$ for each $\gamma\in\Dc$ and $d_k(\phi)=1$ for each $k$.

For each pair  $\ell\geq k>m\geq -\ell$ and $x\in \Uc$, we denote by   $r_{km}(x)$ the rank of submatrix of $x$ with the row system $I_k$ and the column systems  $ I_m$. If $x=x_\Dc$, then we simplify the notation  $r_{km}(\Dc)$. 

\Prop\Num\label{Gequi}. 1) In each  $\Gb$-orbit in $\ux$, there exists an element of the form $x_{\Dc,\phi}$  for some basic subset  $\Dc$ and special $\phi$.\\
2) Two elements $x_{\Dc,\phi}$ and $  x_{\Dc',\phi'}$ with special $\phi$ and $\phi'$ belong to a common $\Gb$-orbit if and only if  $r_{km}(\Dc)=r_{km}(\Dc')$
for any $\ell\geq k>m\geq -\ell$,   and $d_k(\phi)=d_k(\phi')$ for any  $\ell\geq k \geq 1$. 
We say that these basic pairs  are  $\Gb$-equivalent.
\\
\Proof.  \textit{Item 1.} Extend the radical subgroup $\Ub$ of the parabolic subgroup $\Gb$ to the unitriangular subgroup $\tUb=\mathrm{UT}(N,\Fq)$.
 For any $x\in \ux$  the orbit $\Gb\centerdot x$ is a union of  of $\tUb$-orbits. Each $\tUb\centerdot y$ coincides with the intersection 
of  $\tUb y\tUb$ with $\ux$ (see \cite[Theorem 6.10]{SA}). It follows that each $\tUb$-orbit contains  an element of $x_{\Dc,\phi}$ for some basic subset $\Dc$ and $\phi:\Dc\to \Fq^*$ (see \cite[Prop. 3.4]{AN-2}). So, each $\Gb$-orbit contains some  $x_{\Dc,\phi}$.

To prove 1) it remains to show that each  $x_{\Dc,\phi}$ is equivalent to  the element with special $\phi$. Replacing the element  $x_{\Dc,\phi}$ by  $h\centerdot x_{\Dc,\phi}$ with a suitable diagonal matrix $h\in\GL(N,\Fq)$ one may obtain the element  of (\ref{dphi}) with $\phi(\gamma)=1$ for any $\gamma=(i,j)$,~ $j\ne -i$, and  $\phi(i,-i)\in\{1,\delta\}$.

 Finally, let us show that we may choose $\phi$ satisfying  3) of Definition \ref{special}. 
 Suppose that $\Dp$ is the  set of positive roots of type $C_n$.
 
  Denote $X=x_{\Dc,\phi}$. For any $\ell\geq k> m \geq -\ell$, let $X_{k,m}$ be its $I_k\times I_m$-submatrix.  
  Choose  $\ell\geq k\geq 1$.  The $k$th block-row consists of $n_k$ rows as $|I_k|=n_k$. 
  To simplify notations suppose that the row numbers of the $\Dc\cup \Dc'$  belonging to the  blocks $I_k\times I_m$,~ $m\ne -k$,  are the last $n_k''$ row numbers in the segment $I_k$. Take $n_k'=n_k-n_k''$.
  
  Then the submatrix $X_{k,-k}$ has the form
  \begin{equation}\label{kmk}
  X_{k,-k}=\left(\begin{array}{cc}0&Z\\0&0\end{array}\right)
  \end{equation} with some $n_k'\times n_k'$ submatrix $Z$.
  
 Choose non-degenerate  $n_k\times  n_k $-matrix $A_k$ of the form
  $
  \left(\begin{array}{cc}B_k&0\\0&E\end{array}\right)
  $
  with some $n_k'\times n_k'$ submatrix $B_k$ and the identity $n_k''\times n_k''$ submatrix $E$. 
 Consider the  matrix  $R=\diag(A_\ell,\ldots,A_k,\ldots,A_{-\ell})$ of $\Lb$, where each $A_m$, ~$m\ne k$, is the identity matrix, and $A_k$ as above. 
 Then $(R\centerdot X)_{ab}=X_{ab}$ for all $(a,b)\ne (k,-k)$, and 
 the $(k,-k)$-block of  $R\centerdot X$  has the form 
  \begin{equation}\label{Xk}
 \left(\begin{array}{cc}0&B_kZB_k^\dag\\0&0\end{array}\right)
 \end{equation}  
According  the classification of quadratic forms over the finite field $\Fq$ of odd characteristic any quadratic form is equivalent to  the form $x_1^2+\ldots +x_{r-1}^2+x_r^2$ or $x_1^2+\ldots +x_{r-1}^2+\delta x_r^2$, where $r$ is its rank
(see \cite{LN}).
 
Choosing the appropriate non-degenerate matrix  $B_k$  we obtain the matrix $R\centerdot X$ of the from (\ref{dphi}) that obeys condition 3) of Definition \ref{special} in $(k,-k)$ block. This  proves statement 1).\\
\textit{Item 2.1} Let  $X=x_{\Dc,\phi}$, ~ $X'=x_{\Dc',\phi'}$,  and $r_{km}=r_{km}(\Dc)$, ~$r'_{km}=r'_{km}(\Dc')$,~ $d_{k}=d_{k}(\phi)$, ~ 
$d'_{k}=d_{k}(\phi')$.

Suppose that $X$ and $X'$ belong to a common $\Gb$-orbit. Let us prove that $r_{km}=r'_{km}$ for all  $\ell\geq k>m\geq -\ell$, and $d_k=d'_k$  for for all  $\ell\geq k\geq 1$. 
Observe that  for any $x\in \ux$, and  any $\ell\geq k>m\geq -\ell$ the rank of submatrix with the rows and columns $I_k\cup\ldots\cup I_m$   are constant on the orbit $\Gb\centerdot x$. This implies that as $X$ and $X'$ belong to a common $\Gb$-orbit, then $r_{km}=r'_{km}$ for any $\ell\geq k>m\geq -\ell$.

In its turn, the  equality $r_{km}=r'_{km}$ implies that there exist the block-diagonal  matrix  $h=\diag(\sigma_\ell, \ldots,\sigma_{\-\ell})$  with permutations $\sigma_k$, ~$\ell\geq k\geq -\ell$, such that $h\centerdot X'$ has  the same form as $X'$ with $\Dc'=\Dc$.  It remains to show  $d_k=d_k'$ under condition $\Dc'=\Dc$. 

Suppose that $(k,-k)$ blocks in $X$ and $X'$ have the form (\ref{kmk}) with submatrices $Z$ and $Z'$ respectively.
If $X'=g\centerdot X$ then as in Item 1 its $k$th diagonal block has the form
$
\left(\begin{array}{cc}B_k&0\\C&E\end{array}\right)
$ 
for non-degenerate $n_k'\times n_k'$ matrix $B_k$ and some $n_k''\times n_k'$ matrix $C$. Then $Z'=B_k Z B_k^\dag$. This implies $d_k=d_k'$.\\
\textit{Item 2.2} Suppose that $r_{km}=r'_{km}$ and  $d_k=d_k'$ for any $\ell\geq k>m\geq -\ell$, then as in Item 2.1 they coincide under  the action of block-diagonal permutation matrix. $\Box$

In the case of type $A$, the similar statement   follows from the proof of \cite[Theorem 4.2]{Thiem}.

We attach to each basic subset   $\Dc$ and   map $\phi: \Dc\to\Fq$ the element in  $\ux^*$ of the form  $$\la_{\Dc,\phi}=\sum_{\gamma\in\Dc}\phi(\gamma)\Ec_\gamma^*,$$
where  $\{\Ec_\gamma^*\}$ is a dual basis for the basis  $\{\Ec_\gamma\}$ in $\ux$.

The following  statement is analogical to Proposition  \ref{Gequi}.\\
\Prop \Num\label{Gstarequi}.  1) In each  $\Gb$-orbit in $\ux^*$, there exists an element of the form $\la_{\Dc,\phi}$  for some basic subset  $\Dc$ and special $\phi$.\\
2) Two elements $\la_{\Dc,\phi}$ and $  \la_{\Dc',\phi'}$ with special $\phi$ and $\phi'$ belong to a common $\Gb$-orbit if and only if  $r_{km}(\Dc)=r_{km}(\Dc')$
for any $\ell\geq k>m\geq -\ell$,   and $d_k(\phi)=d_k(\phi')$ for any  $\ell\geq k \geq 1$. 

For two decompositions  $\Ic_1$ and $\Ic_2$ of $I$ into union of consecutive segments, we say that $\Ic_1$ \textit{is finer }than $\Ic_2$ (respectively, of $\Ic_2$ \textit{is coarser} than $\Ic_1$) if  each segment of $\Ic_1$ is contained in one of $\Ic_2$ segments.  
By a basic subset  $\Dc$, we define a  new  decomposition   $\Ic^\Dc$ of the segment $I$ into union of consecutive segments such that:\\
1) it is coarser than  $\Ic$, \\
2) it is symmetric about  zero,\\
3)~ for each  $(i,j)\in \Dc$,~ $i\in I_k$,~ $j\in I_m$, the segments  $I_k$ and $I_m$ belong to a common segment in $\Ic^\Dc$,\\
4) it is the finest decomposition  that obeys the above conditions.
 
{\bf Example}. The Borel subgroup   $G$  in $B_2$ is defined by the decomposition  $\Ic=\{\{2\}, \{1\}, \{0\},\{-1\},\{-2\}\}$ of the segment  $I=\{2,1,0,-1,-2\}$.  If $\Dc= \{(2,1)\}$, then $\Ic^\Dc=\{\{2,1\}, \{0\}, \{-1,-2\}$. If  $\gamma=(2,-1)$, then the decompo\-si\-tion  $\Ic^\Dc$ consists of the only segment  $I$.

The decomposition $\Ic^\Dc$  defines the parabolic subgroup  $\Gb^\Dc\supseteq \Gb$ in $\GL(n,\Fq)$ and the parabolic subgroup  $G^\Dc\supseteq G$  in the orthogonal (symplectic) group.  The subgroup  $G^\Dc$ is a semidirect product  $G^\Dc=L^\Dc U^\Dc$, where $U^\Dc$ is the unipotent radical in  $G^\Dc$. Since   $G^\Dc\supseteq G$, we have $U^\Dc\subseteq U$. 
The subgroup  $G\cap L^\Dc$ is a parabolic subgroup in  $L^\Dc$ that is a semidirect product of the Levi subgroup  $L$ and the unipotent radical  $V^\Dc=U\cap L^\Dc$.
The Springer map  $f:\Ub\to \Uc$ bijectively maps  $U^\Dc$ onto
$\Gb$-invariant ideal  $\ux^\Dc$ in $\ux$. Analogically, $f$ bijectively maps  $V^\Dc$ on to some Lie subalgebra  $\vx^\Dc$ in $\ux$.
Observe that  $x_{\Dc,\phi}\in \vx^\Dc$ and $\la_{\Dc,\phi}(\ux^\Dc)=0$. Since the ideal  $\ux^\Dc$ is invariant with respect to the  $\Gb$-action on  $\ux$,  the linear forms from  $\Gb\centerdot\la_{\Dc,\phi}$ annihilate on  $\ux^\Dc$. 

It follows from the Proposition  \ref{uchi} that for any linear form   $\la\in\ux^*$ that is equal to zero on  $\ux^\Dc$ the formula
\begin{equation}\label{guchi}
\zeta_{\Gb\centerdot\la}(u)  = \sum_{\mu\in \Gb\centerdot\la} \eps^{\mu(f(u))}
\end{equation} 
defines a character of  subgroup  $V^\Dc$.  Since $V^\Dc = U/U^\Dc$, the formula  (\ref{guchi}) defines also a character of subgroup  $U$.

Consider the subgroup $L_\Dc$ of all matrices  $\mathrm{diag}(A_\ell,\ldots,A_0,\ldots, A_{-\ell})$ from $L$ such that each its submatrix  $\mathrm(A_k,\ldots,A_m)$ is scalar if  $k\ne m$ and  $I_k\sqcup\ldots \sqcup I_m$ is a segment in decomposition  $\Ic^\Dc$.\\
{\bf Example}. Let $G$ be the Borel subgroup  in  $B_2$, and $\Dc$ consists of a unique root  $\gamma$. If $\gamma=(2,1)$, then  $L_\Dc=\{\diag(a,a,\pm 1,a^{-1},a^{-1})\}$.  If $\gamma=(2,-1)$, then  $L_\Dc=\{\pm E\}$, where $E$ is the identity matrix.\\
\Remark. Easy to see that   $L_\Dc$ coincides with the subgroup  
$L_{\Gb\centerdot\la_{\Dc,\phi}} $ consisting  of all  $h\in L$ such that   $\Ad^*_h(\mu)=\mu$ for all  $ \mu\in\Gb\centerdot\la_{\Dc,\phi}$. In its turn, the last subgroup coincides with the subgroup $L_{\Gb\La_{\Dc,\phi}\Gb} $ consisting of all $h\in L$ such that  $\Ad^*_h(M)=M$ for all  $ M\in\Gb\La_{\Dc,\phi}\Gb$. This implies that  $L_\Dc$ is a subgroup in  $L_{\Ub \La_{\Dc,\phi}\Ub}$.

The formula $$\chi_{\theta,\Dc,\phi}^\circ(g_0)=\theta(r_0) \zeta_{\Gb\centerdot\la_{\Dc,\phi}}(u),$$ where $g_0=r_0u$, ~$r_0\in L_D$, ~$u\in U$,
defines a character of the subgroup  $L_\Dc U$.

Consider the set  $$\Ac =\{(\theta,\Dc,\phi)\}, $$ where $(\Dc,\phi)$ runs through the set of   representatives of $\Gb$-equivalency classes of basic pairs, and $\theta$  runs through the set of all irreducible characters of the group  $L_{\Dc}$.

We attach to each   $\al=(\theta,\Dc.\phi)$ the character  $\chi_\al$ induced from the character $\chi_{\theta,\Dc,\phi}^\circ$ of subgroup  $L_\Dc U$ to $G=LU$. 

The subgroup $L_{\Dc}$  is a normal subgroup in $L$. Moreover, for each character   $\theta$ of the subgroup  $L_\Dc$, the following property is fulfilled: $\theta(\rho r \rho^{-1})=\theta(r)$ for any  $r, \rho\in L_\Dc$.  The character  $\chi_\al $ is calculated by the formula 

\begin{equation}\label{chialal}
\chi_\al(g)=\frac{|L|}{|L_\Dc|}\cdot \dot{\theta}(r)\cdot \sum_{\mu\in \Gb\centerdot\la_{\Dc,\phi}} \eps^{\mu(f(u))},
\end{equation}
where $g=ru$,~ $r\in L$,~ $u\in U$, and $\dot{\theta}$ is the function on  $L$ that is equal to  $\theta$ on $L_\Dc$ and is equal to zero outside of  $L_\Dc$.

Consider the set  $$\Bc=\{(h,\Dc,\phi)\},$$
where $(\Dc,\phi)$ runs through the set of   representatives of $\Gb$-equivalency classes of basic pairs, and 
$h$ runs through the set of representatives of conjugacy classes in $L_\Dc$.

By the element  $h=\diag(A_\ell,\ldots ,A_0,\ldots, A_{-\ell})\in L$, we define  a new decom\-po\-si\-tion  $\Ic_h$  of  the segment $I$ into union of consecutive segments such that:\\
1) it is coarser than  $\Ic$, \\
2) it is symmetric about zero,\\
3)  the segments  $I_k$ and $I_m$, where $k> m$, belong to a common segment in  $\Ic_h$ if $\diag(A_k,\ldots,A_m)$ is a scalar submatrix of  $h$,\\
4) it is the coarsest  decomposition  that obeys the above conditions.

The decomposition $\Ic_h$ defines the parabolic subgroup $\Gb_h$ in $\GL(m,\Fq)$, which  contains  $\Gb$, and  the parabolic subgroup   $G_h=L_hU_h$ in the orthogonal (symplectic) group,  which contains  $G$. 
The subgroup  $U$ is a semidirect product  $U=V_hU_h$, where $V_h=L\cap U$. 
The Lie algebra  $\vx_h$ of the subgroup $V_h$ is spanned by $\{\Ec_\gamma:~  h\in L_\gamma\}$.
Respectively, $\ux_h$ is spanned by  $\{\Ec_\gamma:~  h\notin L_\gamma\}$.

We attach to each element  $\beta=(h,\Dc)\in \Bc$ the subset  
$$ K_\beta = \mathrm{Cl}(h)\cdot K_{\Dc,\phi}\cdot U_h,$$
where $\mathrm{Cl}(h)$ is the conjugacy class for $h$ in $L_\Dc$, 
and $K_\Dc = 1+ \Gb\centerdot x_{\Dc,\phi}$. \\
\Theorem\Num\label{superchar}.  The system of characters $\{\chi_\al:~\al\in\Ac\}$ and the partition of the parabolic subgroup  $G$ into the classes  $\{K_\beta:~\beta\in \Bc\}$ give rise to a supercharacter theory.  \\
\Proof. The characters from  $\{\chi_\al\}$  attached to different  $\Gb$-orbits in  $\ux^*$ are orthogonal since their restrictions on the normal subgroup  $U$ are orthogonal. It follows from the formula  (\ref{chialal}) that the characters attached to a common  $\Gb$-orbit and different irreducible characters of the subgroup  $L_\Dc$ are orthogonal.

The value of the character  $\chi_\al$, ~ $\al=(\theta, \Dc,\phi)$, on the class  $K_\beta$,~
$\beta=(h,\Dc',\phi')$, is calculated by the formula 
$$\chi_\al(K_\beta) = \dot{\theta}(h) \zeta_{\Gb\centerdot\la_{\Dc,\phi}}(K_{\Dc',\phi'}) =\mathrm{const}.$$
Finally,  $|\Ac|=|\Bc|$ because for each $\Dc$ the number of conjugacy classes  in $L_\Dc$ equals  to the number of irreducible representations of $L_\Dc$.  $\Box$
\\
\Cor\Num. 1) Each irreducible character of $G$ is a constituent of a unique supercharacter   from  $\{\chi_\al:~\al\in\Ac\}$.
2) Each conjugacy class of $G$ is contained in a unique class from    $\{K_\beta:~\beta\in \Bc\}$.

\end{document}